\documentclass[a4,amstex,11pt]{article}
\usepackage{amsmath,wrapfig}
\usepackage{amssymb,bm,epic}
\usepackage[mathscr]{eucal}
\usepackage{amsthm, amsfonts, latexsym,enumerate}
\usepackage{graphicx}
\usepackage{setspace}
\makeatletter
\def\fixedfigure{\def\@captype{figure}}
\def\fixedtable{\def\@captype{table}}
\numberwithin{equation}{section}
\theoremstyle{definition}

\newtheorem{definition}{Definition}[section]
\newtheorem{theorem}[definition]{Theorem}
\newtheorem{proposition}[definition]{Proposition}
\newtheorem{lemma}[definition]{Lemma}
\newtheorem{corollary}[definition]{Corollary}

\newtheorem{example}[definition]{Example}

\newtheorem{algorithm}[definition]{Algorithm}
\theoremstyle{break}

\title{On max-plus two-sided linear systems whose solution sets are min-plus linear
\thanks{This work is supported by JSPS KAKENHI Grant No.~22K13964.}
}

\author{
Yasutaka Ooga\footnote{Science of Environment and Mathematical Modeling, Doshisha University, Kyoto, Japan.}
\and
Yuki Nishida\footnote{Department of Information and Computer Technology, Tokyo University of Science, Tokyo, Japan.
Email: ynishida.cyjc1901@gmail.com
}
\and 
Yoshihide Watanabe\footnote{Department of Mathematical Sciences, Doshisha University, Kyoto, Japan.}
}
\date{}
\begin{document}
\maketitle


\begin{abstract}
The max-plus algebra $\mathbb{R}\cup \{-\infty \}$ is defined in terms of a combination of the following two operations: 
addition, $a \oplus b := \max(a,b)$, and multiplication, $a \otimes b := a + b$.
In this study, we propose a new method to characterize the set of all solutions 
of a max-plus two-sided linear system $A \otimes \bm{x} = B \otimes \bm{x}$.
We demonstrate that the minimum ``min-plus'' linear subspace containing the ``max-plus'' solution space can be computed by applying the alternating method algorithm, which is a well-known method to compute single solutions of two-sided systems.
Further, we derive a sufficient condition for the ``min-plus'' and ``max-plus'' subspaces to be identical.
The computational complexity of the method presented in this study is pseudo-polynomial.
\end{abstract}



{\it 2020MSC}: 15A80, 15A06, 15A39

\section{Introduction}
\label{sec1}
The max-plus algebra, $\mathbb{R}\cup \{-\infty \}$, is defined in terms of the following two operations: 
addition, $a \oplus b := \max(a,b)$, and multiplication, $a \otimes b := a + b$.
The max-plus algebra originated in steelworks~\cite{Green1960, Green1962}.
It has a wide range of applications in various fields of science and engineering, such as control theory and
scheduling in railway systems~\cite{Baccelli1992, Heidergott2005}. 
\par
Linear systems over the max-plus algebra have been a popular topic of research for a long time.
One-sided systems $A \otimes \bm{x} = \bm{b}$ were solved by the combinatorial way in the 1960s~\cite{Green1960}.
They were later solved algebraically by introducing a type of dual semiring with two operators: 
$a \oplus' b := \min(a,b)$ and $a \otimes' b := a + b$~\cite{Green1976}.
However, two-sided systems $A \otimes \bm{x} = B \otimes \bm{x}$ are difficult to solve
and remain an active topic of research from both theoretical and computational perspectives.
A two-sided system can be reduced to a separated system $A \otimes \bm{x} = B \otimes \bm{y}$, which can be solved using the alternating method~\cite{Green2003}.
Based on the algebraic way to solve one-sided systems,
this algorithm alternately applies max-plus and min-plus matrix multiplications.
The alternating method is a pseudo-polynomial time algorithm if either $A$ or $B$ is an integer matrix.
Currently, there is no known polynomial-time algorithm that can solve general two-sided max-plus linear systems.
The existence of a solution to a two-sided system is equivalent to a mean payoff game~\cite{Akian2012}.
Several algorithms have been developed to investigate the general case and its complexities~\cite{Bezem2008, Butkovic2006, Grigoriev2015}.
Tropical linear systems, e.g., $A \otimes \bm{x}$, have also been considered in the literature,
in which the solution is a vector such that the maximum $\bigoplus_{j} a_{ij} \otimes x_{j}$ is attained at least twice in each row $i$.
This type of solution is obtained based on the theory of tropical geometry over fields with valuation~\cite{Joswig2021, Maclagan2015}.
The solutions of $A \otimes \bm{x} = B \otimes \bm{x}$ are closely related to those of a tropical linear system $(A \oplus B) \otimes \bm{x}$.
Further, tropical linear systems can be reduced to max-plus two-sided linear systems.
Several algorithms have been proposed to investigate tropical linear systems~\cite{Davydow2013, Davydow2017, Grigoriev2013}.
In certain cases, either max-plus two-sided systems or tropical linear systems can be solved more efficiently.
When the matrix $A$ is of size $n\times (n+1)$,
tropical Cramer's rule can be used to derive a tropical solution in $O(n^{3})$ time~\cite{Akian2014,Gebert2005}.
Cramer's rule has been extended to overdetermined cases~\cite{Davydow2018}.
When $A$ and $B$ are square matrices, two-sided max-plus linear systems can be solved under the assumption of idempotency of matrices~\cite{Butkovic1985}, strong T systems~\cite{Aminu2010}, 
or minimally active or essential systems~\cite{Jones2019}.
Linear systems can also be solved based on the symmetrized max-plus algebra by augmenting ``sign-negative'' elements and 
``balanced elements,” which behave as zeros~\cite{Plus1990}. 
Then, each term can be moved to the other side of the equality, yielding 
$(A \ominus B) \otimes \bm{x} = \bm{\varepsilon}$.
In particular, this method is useful for solving the system $A \otimes \bm{x} \oplus \bm{c} = B \otimes \bm{x} \oplus \bm{d}$ for square matrices $A$ and $B$ because the inverse of $(A \ominus B)$ can be considered in terms of  the symmetrized max-plus algebra.
The concept of the symmetrized max-plus algebra has been extended to supertropical algebra and 
several techniques have been developed to investigate it~\cite{Izhakian2011}.
\par
The determination of the set of all solutions of max-plus two-sided linear systems is also an important problem.
One method to identify all solutions using the usual linear inequalities was presented in~\cite{Lorenzo2011}.
Because the solution set is closed under addition and scalar multiplication in max-plus arithmetic,
it is sufficient to determine a generating set of solutions as a max-plus subspace.
A fundamental elimination method for deriving a generating set of the solution space was proposed in~\cite{Butkovic1984}.
Although the elimination step is exhaustive and difficult to execute for large matrices, it ensures that the solution set is a finitely generated max-plus subspace.
This idea is applied to derive efficient algorithms with a small number of inequalities~\cite{Sergeev2011, Wagneur2009}
and the sparsification method for solving max-plus linear inequalities~\cite{Krivulin2020}.
Via analogy with polyhedral cone theory,
the tropical double description method was proposed to compute a generating set of the solution space~\cite{Allamigeon2013}.
Additionally, tropical Cramer's rule can be used to solve this problem~\cite{Nishida2020}.
The execution of such algorithms requires exponential time in the worst-case scenario.
The known upper bound for the cardinality of the basis of the solution space can be expressed in terms of binomial coefficients
in the sizes of matrices~\cite{Allamigeon2011}.
The method presented in~\cite{Goncalves2012} efficiently finds a large number of solutions based on a given solution.
\par
The purpose of this study is to characterize the solution set of 
a max-plus two-sided linear system $A \otimes \bm{x} = B \otimes \bm{x}$ from the min-plus linear algebraic perspective.
The primary tool used is the alternating method ~\cite{Green2003}.
As described above, the alternating method computes a single solution for each two-sided system.
The obtained solution depends on the initial vectors.
Hence, every solution may be obtained via the alternating method if the initial vector is appropriately selected.
In the present study, we consider $m$ initial vectors corresponding to the $m$ rows of the system.
Then, we obtain $m$ solutions using the alternating method 
and demonstrate that all the solutions of the system are contained in the min-plus linear subspace generated by these $m$ vectors.
In other words, the min-plus closure of the solution set can be computed.
In particular, if the solution set  is also closed under the ``min'' operation for vectors, 
the solution set is then completely determined by applying the alternating method $m$ times.
This type of subset plays an important role in combinatorial optimization.
A subset that is both max-plus and min-plus linear is called an L-convex set~\cite{Murota1998}.
In this study, we also derive a sufficient condition for the solution set of a two-sided system to be min-plus linear.
Although the min-plus linearity of a subset is a global condition,
we demonstrate that local min-plus convexity around the $m$ vectors computed using the alternating method is a sufficient condition for it.
This implies that the complexity of computing the min-plus linear closure of the solution set as well as that of verifying the min-plus linearity of the solution set are both pseudo-polynomial time for two-sided systems defined by integer matrices.
\par
The remainder of this paper is organized as follows.
In Section 2, we introduce basic definitions of max-plus and min-plus linear algebra.
In section 3, we summarize the results on max-plus two-sided linear systems and the solutions obtained using the alternating method.
Stable solutions with respect to the algorithm play an important role in the discussion.
In Section 4, we first restrict the set of solutions to those with finite entries to ensure that the solution set is bounded in the max-plus projective space.
Then, we prove that the vectors computed via the alternating method generate the minimum min-plus linear subspace containing all solutions.
In Section 5, we present a criterion for the solution set to be min-plus linear.
First, we observe that min-plus linearity of a max-plus subspace is characterized in terms of local min-plus convexity.
Then, we describe the cell where local min-plus convexity of a solution is violated.
Finally, we demonstrate that verification of local min-plus convexity at the vectors computed via the alternating method is sufficient. 
\section{Max-plus algebra}
\label{sec2}
\subsection{Max-plus and min-plus algebras}
Let $\mathbb{R} _{\max} = \mathbb{R} \cup \{\varepsilon\}$ be the set of real numbers $\mathbb{R}$ with an extra element
 $\varepsilon:=-\infty$.
We define two operations---addition, $\oplus$, and multiplication, $\otimes$---on $\mathbb{R} _{\max}$ as follows:
\begin{align*}
a \oplus b=\max(a,b), \quad 
a \otimes b=a+b, \quad 
a,b \in\mathbb{R} _{\max}.
\end{align*}
Then, ($\mathbb{R}_{\max},\oplus,\otimes$) is a commutative semiring called the max-plus algebra. 
Here, $\varepsilon$ is the additive identity and $e := 0$ is the multiplicative identity.
\par
Similarly, the min-plus algebra is denoted by $\mathbb{R}_{\min} = \mathbb{R} \cup \{\varepsilon'\}$, where $\varepsilon' := \infty$ and the following two operations hold:
\begin{align*}
a \oplus' b=\min(a,b), \quad
a \otimes' b=a+b, \quad 
a,b \in\mathbb{R} _{\min}.
\end{align*}
We extend the operations $\oplus, \otimes, \oplus'$, and $\otimes'$ to the set 
$\mathbb{R}_{\max,\min} := \mathbb{R} \cup \{\varepsilon,\varepsilon'\}$ as follows:
\begin{align*}
	\varepsilon \oplus \varepsilon' &= \varepsilon' \oplus \varepsilon = \varepsilon', \quad
	\varepsilon \oplus' \varepsilon' = \varepsilon' \oplus' \varepsilon = \varepsilon, \\
	\varepsilon \otimes \varepsilon' &= \varepsilon' \otimes \varepsilon = \varepsilon, \quad
	\varepsilon \otimes' \varepsilon' = \varepsilon' \otimes' \varepsilon = \varepsilon'.
\end{align*}
For further details regarding max-plus and min-plus algebras,
please refer to~\cite{Baccelli1992, Butkovic2010, Green1979, Gondran2010, Heidergott2005, Joswig2021, Maclagan2015}.
\subsection{Max-plus linear algebra}
Let $\mathbb{R}_{\max}^n$ and $\mathbb{R}_{\max}^{m\times n}$ be sets of $n$-dimensional max-plus column vectors 
and $m\times n$ max-plus matrices, respectively. 
The operations $\oplus$ and $\otimes$, are extended to max-plus vectors and matrices following conventional linear algebra. 
Let $\bm{\varepsilon}$ and $\mathcal{E}$ denote the max-plus zero vector and the zero matrix, respectively.
Further, let $E_n$ denote the max-plus unit matrix of order $n$.
\par
A subset of $\mathbb{R}_{\max}^n$ is called a max-plus subspace if it is closed under (max-plus) addition and scalar multiplication.
For a subset $S\subset \mathbb{R}_{\max}^n$, we can easily verify that the set
\begin{align*}
	\langle S \rangle_{\max} = \left\{ \bigoplus^{r}_{i=1} a_i \otimes \bm{x}_i  \biggm| 
	r \in \mathbb{Z}_{>0}, \ \bm{x}_i\in S, \ a_i\in\mathbb{R}_{\max} \text{ for } i=1,2,\dots,r \right\}
\end{align*}
is a max-plus subspace of $\mathbb{R}_{\max}^{n}$.
Subsequently, $\langle S \rangle_{\max}$ is called the max-plus subspace of $\mathbb{R}^{n}_{\max}$ generated by $S$. 
In particular, the max-plus subspace generated by a finite set $\{ \bm{x}_{1}, \bm{x}_{2}, \dots, \bm{x}_{r} \}$ is denoted by $\langle \bm{x}_{1}, \bm{x}_{2}, \dots, \bm{x}_{r} \rangle_{\max}$.
The minimum generating set of a max-plus subspace is called its basis.
In the max-plus algebra, any finitely generated subspace has a basis that is unique up to scalar multiplication~\cite{Butkovic2007}.
\par
In the present study, we say that the max-plus subspace $V$ is projectively bounded
 if $V \subset \mathbb{R}^n \cup \{\bm{\varepsilon}\}$.
\par
The definitions and notions describe in this subsection extend to analogous ones corresponding to the min-plus algebra.
Thus, $\mathbb{R}_{\min}^n$ denotes the set of min-plus vectors, $\bm{\varepsilon}'$ denotes the min-plus zero vector, etc.
A vector or matrix is considered to be finite if it does not contain $\varepsilon$ or $\varepsilon'$ as an element.

\section{Max-plus linear systems}
\label{3}
\subsection{Alternating method for solving two-sided systems}
The purpose of this study is to investigate all solutions of max-plus homogeneous linear systems
\begin{align}
	A \otimes \bm{x} = B \otimes \bm{x} \label{hom-sys}
\end{align}
with $A, B \in \mathbb{R}_{\max}^{m \times n}$.
Henceforth, we assume that the coefficient matrices, $A$ and $B$, are doubly $\mathbb{R}$-astic, 
i.e., all rows and columns contain at least one real number.
Let $S(A,B) \subset \mathbb{R}_{\max}^{n}$ be the set of all solutions of \eqref{hom-sys}.
Then, $S(A,B)$ is a finitely generated max-plus subspace~\cite{Butkovic1984}.
Although various algorithms have been developed to solve max-plus linear systems,
we focus on an iterative algorithm called the alternating method.
\par
First, consider a max-plus linear system with separate variables
\begin{align}
	A\otimes\bm{x}=B\otimes\bm{y}, \label{sep-sys}
\end{align}
where $A \in \mathbb{R}_{\max}^{m\times n}, B \in \mathbb{R}_{\max}^{m\times k}$. 
Algorithm \ref{A-M} solves the system \eqref{sep-sys}.
We introduce the following notations.
The set $\{1,2,\dots,n\}$ is denoted by $[n]$.
In a sequence of vectors $\{\bm{x}(r)\}_{r=0,1,\dots}$, the $i$th entry of $\bm{x}(r)$ is denoted by $x_{i}(r)$.
Given $A = (a_{ij}) \in \mathbb{R}_{\max}^{m \times n}$, the matrix $-A^{T} = (-a_{ji}) \in \mathbb{R}_{\min}^{n \times m}$
 is sometimes used as its pseudo-inverse matrix. 
\begin{algorithm}[Alternating Method \cite{Green2003}]\label{A-M}
\begin{enumerate}
\renewcommand{\labelenumi}{(\theenumi)}
\item Take any vector $\bm{x}(0) \in \mathbb{R}^n$ and set $r:=0$.
\item We compute
\begin{align*}
\bm{y}(r)  := -B^T \otimes' (A \otimes \bm{x}(r)), \quad
\bm{x}(r + 1)  :=-A^T \otimes' (B \otimes \bm{y}(r)).
\end{align*}
\begin{enumerate}
\renewcommand{\labelenumi}{(\alph{\enumi})}
\item If $x_i (r + 1) < x_i (0)$ for $i \in [n]$, then there is no finite solution of \eqref{sep-sys}.
\item If $\bm{x}(r+1) = \bm{x}(r)$, then $(\bm{x}(r),\bm{y}(r))$ is a solution of \eqref{sep-sys} in $\mathbb{R}^n$.
\item If neither (a) nor (b) holds, then set $r := r+1$ and repeat (2).
\end{enumerate}
\end{enumerate}
\end{algorithm}
The alternating method is a pseudo-polynomial algorithm.
The following results demonstrate the computational complexity of the algorithm.
\begin{proposition}[\cite{Green2003}] \label{A-Mcomplex}
If $A\in\mathbb{Z}^{m\times n},B\in(\mathbb{Z}\cup\{\varepsilon\})^{m\times k}$ and $\bm{x}(0) \in\mathbb{Z}^n$, 
then the alternating method terminates after finitely many steps.
Further, if  $|x_{i}(0)|$ is bounded by $K:=\max_{i,j} |a_{ij}|$, the computational complexity is $O(mn(n+k)K)$.
\end{proposition}
Now, we consider the max-plus linear system \eqref{hom-sys}.
This system is equivalent to the following system with separate variables:
\begin{align}
	\begin{pmatrix}
	A\\
	B
	\end{pmatrix}
	\otimes\bm{x}=
	\begin{pmatrix}
	E_m \\
	E_m
	\end{pmatrix}
	\otimes\bm{y} \label{hom-sep}.
\end{align}
Hence, we can obtain a finite solution of \eqref{hom-sys} by applying the alternating method to \eqref{hom-sep}.
When $A = (a_{ij}), B = (b_{ij}) \in \mathbb{Z}^{m \times n}$,
the computational complexity of the alternating method for~\eqref{hom-sys} is $O(mn(m+n)K)$,
where $K=\max(\max_{i,j} |a_{ij}|, \max_{i,j} |b_{ij}|)$.
One iteration step of the Algorithm \ref{A-M} for \eqref{hom-sep} proceeds as follows:
\begin{align*}
	\varphi_{0}(\bm{x}) :=& 
	\begin{pmatrix} -A^{T} & -B^{T} \end{pmatrix} \otimes' \left( \begin{pmatrix} E_{m} \\ E_{m} \end{pmatrix}
	\otimes \left( \begin{pmatrix} -E_{m}^{T} & -E_{m}^{T} \end{pmatrix} \otimes' \left( \begin{pmatrix} A \\ B \end{pmatrix}
	\otimes \bm{x} \right) \right) \right) \\
	=& (-A^{T} \oplus' -B^{T}) \otimes' ( (A \otimes \bm{x}) \oplus' (B \otimes \bm{x}) ). 
\end{align*}
\begin{example} \label{exm1}
	Consider the linear system \eqref{hom-sys} for
	\begin{align*}
		A = \begin{pmatrix} 0  & 1 & -1 \\ 0 & -5 & -5 \\ 0 & 4 & 6 \\ 0 & 3 & -2 \end{pmatrix}, \quad
		B = \begin{pmatrix} 0 & -1 & -1 \\ 0 & -4 & -3 \\ -1 & 1 & 6 \\ -1 & 3 & -3 \end{pmatrix}.
	\end{align*} 
	Let $\bm{x}(0) = (0,4,3)^{T}$ be the initial vector of the alternating method.
	Then, we have:
	\begin{align*}
		\bm{x}(1) &= \begin{pmatrix} 0 & 0 & 0 & 0 \\ -1 & 4 & -4 & -3 \\ 1 & 3 & -6 & 2 \end{pmatrix}
			\otimes' \left( \begin{pmatrix} 5 \\ 0 \\ 9 \\ 7 \end{pmatrix} 
			\oplus' \begin{pmatrix} 3 \\ 0 \\ 9 \\ 7 \end{pmatrix} \right)
		= \begin{pmatrix} 0 \\ 2 \\ 3 \end{pmatrix}, \\
		\bm{x}(2) &= \begin{pmatrix} 0 & 0 & 0 & 0 \\ -1 & 4 & -4 & -3 \\ 1 & 3 & -6 & 2 \end{pmatrix}
			\otimes' \left( \begin{pmatrix} 3 \\ 0 \\ 9 \\ 5 \end{pmatrix} 
			\oplus' \begin{pmatrix} 2 \\ 0 \\ 9 \\ 5 \end{pmatrix} \right)
		= \begin{pmatrix} 0 \\ 1 \\ 3 \end{pmatrix}, \\
		\bm{x}(3) &= \begin{pmatrix} 0 & 0 & 0 & 0 \\ -1 & 4 & -4 & -3 \\ 1 & 3 & -6 & 2 \end{pmatrix}
			\otimes' \left( \begin{pmatrix} 2 \\ 0 \\ 9 \\ 4 \end{pmatrix} 
			\oplus' \begin{pmatrix} 2 \\ 0 \\ 9 \\ 4 \end{pmatrix} \right)
		= \begin{pmatrix} 0 \\ 1 \\ 3 \end{pmatrix}.
	\end{align*}
	Hence, $\bm{x} = (0,1,3)^{T}$ is the solution to \eqref{hom-sys}.
\end{example}
\subsection{Stable solutions obtained via the alternating method}
Let us consider the max-plus linear system \eqref{hom-sys}.
By applying Algorithm \ref{A-M} to the max-plus linear system \eqref{hom-sep},
the monotonicity and stability of the iteration are evident.
Let $\bm{x}(r)$ denote the vector in $r$th iteration of Algorithm \ref{A-M}.
\begin{lemma}[\cite{Green2003}] \label{1-1}
	The sequence $\{\bm{x}(r)\}_{r=1,2,\cdots}$ is non-increasing, i.e., $\bm{x}(r+1) \leq \bm{x}(r)$ for $r \geq 1$.
\end{lemma}
\begin{lemma}[\cite{Green2003}] \label{1-2}
	If $\bm{x}$ is a solution of \eqref{hom-sys}, then $\varphi_{0}(\varphi_{0}(\bm{x})) = \varphi_{0}(\bm{x})$.
\end{lemma}
A solution to \eqref{hom-sys} that satisfies $\varphi_{0}(\bm{x}) = \bm{x}$ is called a stable solution.
The solution obtained via the alternating method is stable.
As $\bm{x} \in S(A,B)$ implies $A \otimes \bm{x} = B \otimes \bm{x}$, 
a stable solution also satisfies 
\begin{align}
	\bm{x}=(-A^T)\otimes'(A\otimes\bm{x})\oplus'(-B^T)\otimes'(B\otimes\bm{x}). \label{stable}
\end{align}
The set of all stable solutions is denoted by $\tilde{S}(A,B)$.
\par
Let $A_{i}$ and $B_{i}$ denote the $i$th rows of $A$ and $B$, respectively.
For the system \eqref{hom-sys} and each $i \in [m]$, we define the following set: 
\begin{align*}
M^i(\bm{x})=\{ k\in[n]\ |\ A_i\otimes\bm{x}=a_{ik}\otimes x_k\ \text{or}\ B_i\otimes\bm{x}=b_{ik}\otimes x_k\}.
\end{align*}
Then, stable solutions are characterized by the following proposition:
\begin{proposition} \label{4-1}
	If the solution $\bm{x} \in S(A,B)$ satisfies $\bigcup_{i=1}^m M^i(\bm{x})=[n]$, then $\bm{x}$ is stable.
\end{proposition}
\proof
	We derive the equality $\bm{x}=(-A^T)\otimes'(A\otimes\bm{x})\oplus'(-B^T)\otimes'(B\otimes\bm{x})$.
	The right-hand side of this equality can be expanded as follows:
	\begin{align*}
		&(-A^T)\otimes'(A\otimes\bm{x})\oplus'(-B^T)\otimes'(B\otimes\bm{x}) \\
		=& \bigoplus_{i=1}^m\!\text{\Large ${}'$ } \left((-A_i^T) \otimes' (A_i \otimes \bm{x})
			\oplus' (-B_i^T) \otimes' (B_i \otimes \bm{x}) \right).
	\end{align*}
	For $i \in [m]$, let $k \in [n]$ be an index such that
	$A_{i} \otimes \bm{x} = a_{ik} \otimes x_{k}$ or $B_{i} \otimes \bm{x} = b_{ik} \otimes x_{k}$.
	Then, we have either
	\begin{align*}
		[(-A_i^T) \otimes' (A_i \otimes \bm{x})]_k = -a_{ik} \otimes' (a_{ik} \otimes x_k) = x_k
	\end{align*}
	or
	\begin{align*}
		[(-B_i^T) \otimes' (B_i \otimes \bm{x})]_k = -b_{ik} \otimes' (b_{ik} \otimes x_k) = x_k,
	\end{align*}
	where $[*]_{k}$ denotes the $k$th entry of the vector.
	Therefore, we have
	\begin{align*}
		[(-A_i^T) \otimes' (A_i \otimes \bm{x}) \oplus' (-B_i^T) \otimes' (B_i\otimes\bm{x})]_k \leq x_k.
	\end{align*}
	Hence, the assumption $\bigcup_{i=1}^m M^i(\bm{x})=[n]$ implies
	\begin{align*}
		\bigoplus_{i=1}^m\!\text{\large ${}'$ } (-A_i^T)\otimes'(A_i\otimes\bm{x})\oplus'(-B_i^T)\otimes'(B_i\otimes\bm{x})\leq \bm{x}.
	\end{align*}
	However, we can easily verify that any vector $\bm{x}$ satisfies
	\begin{align*}
		(-A^T) \otimes' (A\otimes\bm{x}) \geq \bm{x}
	\end{align*}
	and that a similar inequality holds for $B$.
	Hence, we have
	\begin{align*}
		\bm{x} \leq (-A^T) \otimes' (A \otimes \bm{x}) \oplus' (-B^T) \otimes' (B \otimes \bm{x}).
	\end{align*}
	This completes the proof of the proposition. 
\endproof

\section{Description of solutions as min-plus linear subspaces}
For any $A, B \in \mathbb{R}_{\max}^{m \times n}$, a finite solution of the max-plus linear system \eqref{hom-sys}
can be determined using the alternating method.
However, determining the solution set $S(A,B)$ is a rather difficult problem.
We demonstrate that the alternating method can be used to resolve this problem effectively, especially when $S(A,B)$ is also min-plus linear.

\subsection{Projectively bounded solution spaces}\label{S-E}
For $\alpha,\beta\in\mathbb{R}$, consider the following matrix:
\begin{align*}
	D_{\alpha,\beta}:=
	\begin{pmatrix}
	\alpha & \beta & \cdots &\cdots & \beta\\
	\beta & \alpha &     &    & \vdots \\
	\vdots &        & \ddots &        & \vdots \\
	\vdots &        &        & \alpha & \beta \\
	\beta & \cdots & \cdots & \beta & \alpha
	\end{pmatrix}.
\end{align*}
Then, the following lemma holds.
\begin{lemma} \label{Dsol}
	For any $\alpha>0$, we have
	\begin{align*}
		S(D_{\alpha,0},D_{\alpha,-1}) = 
			\{\bm{x}\in\mathbb{R}^n\ |\  |x_j-x_k|\leq\alpha \text{ for } j,k\in[n]\} \cup \{\bm{\varepsilon}\}.
	\end{align*}
\end{lemma}
\proof
	For any vector $\bm{x} \in \mathbb{R}^{n}$,
	we have the following equivalence relations:
	\begin{align*}
		& \bm{x}\in S(D_{\alpha,0},D_{\alpha,-1}) \\
		\Longleftrightarrow \quad& \alpha \otimes x_k \oplus \bigoplus_{i \neq k} x_i
			=\alpha \otimes x_k \oplus \bigoplus_{i \neq k}(-1) \otimes x_i \quad ({}^\forall k \in [n])\\
		\Longleftrightarrow \quad& \alpha \otimes x_k \geq \bigoplus_{i \neq k} x_i \quad ({}^\forall k \in [n])\\
		\Longleftrightarrow \quad& \alpha \otimes x_k \geq x_j \quad ({}^\forall j,k \in [n]) \\
		\Longleftrightarrow \quad& \alpha \geq x_j-x_k \quad ({}^\forall j,k \in [n]) \\
		\Longleftrightarrow \quad& |x_j-x_k| \leq \alpha \quad ({}^\forall j,k \in [n]).
	\end{align*}
	We can also verify that all solutions except $\bm{\varepsilon}$ are finite.
\endproof
When we consider the system \eqref{hom-sys},
it is helpful to introduce the matrix $C = (c_{ij}) = -(A \oplus B) \in \mathbb{R}_{\min}^{m \times n}$.
\begin{proposition} \label{Dinc}
	For $A, B\in\mathbb{R}_{\max}^{m\times n}$, 
	let $\alpha$ be any real number larger than $\max\{ c_{ij} - c_{ik} \mid i \in [m],\ j,k \in [n],\ c_{ij}, c_{ik} \in \mathbb{R} \}$.
	Then, the solution set  $S(A,B)$ is projectively bounded if and only if $S(A, B) \subset S(D_{\alpha,0}, D_{\alpha,-1})$.
\end{proposition}
\proof
	``If'' part: By Lemma \ref{Dsol}, $S(D_{\alpha,0}, D_{\alpha,-1})$ is projectively bounded.
	Hence, $S(A,B)$ is projectively bounded if $S(A, B) \subset S(D_{\alpha,0}, D_{\alpha,-1})$.
	\par
	``Only if'' part:
	Suppose $S(A,B)$ is projectively bounded, but $S(A, B) \not\subset S(D_{\alpha,0}, D_{\alpha,-1})$.
	Then, there exists a vector $\bm{x} = (x_{1},x_{2},\dots,x_{n})^{T} \in S(A,B)$ such that $\bm{x} \not\in S(D_{\alpha,0},D_{\alpha,-1})$.
	By Lemma \ref{Dsol}, there exist two indices $j,k \in [n]$, such that
	\begin{align*}
		x_j-x_k > \alpha.
	\end{align*}
	Because $\alpha > \max_{i,j,k} \{ c_{ij} - c_{ik} \}$, we have
	\begin{align*}
		-c_{ij} \otimes x_{j} > -c_{ik} \otimes x_{k}
	\end{align*}
	for all $i \in [m]$.
	Therefore, we have 
	\begin{align*}
		(A_{i} \oplus B_{i}) \otimes \bm{x} > (a_{ik} \oplus b_{ik}) \otimes x_{k}
	\end{align*}
	for all $i \in [m]$.
	Hence, $(x_1,x_2,\dots,x_{k-1},\varepsilon,x_{k+1},\dots,x_n)^T$ is also present in $S(A,B)$.
	This result contradicts the fact that $S(A,B)$ is projectively bounded.
\endproof

Let us consider matrices
\begin{align}
	\tilde{A}_{\alpha} =\begin{pmatrix}
	A\\
	D_{\alpha,0}
	\end{pmatrix}, \quad
	\tilde{B}_{\alpha} =\begin{pmatrix}
	B\\
	D_{\alpha,-1}
	\end{pmatrix}. \label{extmat}
\end{align}
From Proposition \ref{Dinc}, we have $S(A,B) = S(\tilde{A}_{\alpha}, \tilde{B}_{\alpha})$ if $S(A,B)$ is projectively bounded.
Even when $S(A,B)$ is not projectively bounded, for any $\bm{x} \in S(A,B) \cap \mathbb{R}^{n}$, there exists a real number $\alpha > 0$ 
such that $\bm{x} \in S(\tilde{A}_{\alpha}, \tilde{B}_{\alpha})$.
Hence, the set $S(\tilde{A}_{\alpha}, \tilde{B}_{\alpha})$ can be considered to be an approximation of $S(A,B)$ for large $\alpha$.
The system $\tilde{A}_{\alpha} \otimes \bm{x} = \tilde{B}_{\alpha} \otimes \bm{x}$ exhibits good properties when the alternating method is applied to it.
\begin{proposition} \label{stableprop}
	For any $\alpha > 0$, every solution to $\tilde{A}_{\alpha} \otimes \bm{x} = \tilde{B}_{\alpha} \otimes \bm{x}$ is stable.
\end{proposition}
\proof
	Consider any vector, $\bm{x} \in S(\tilde{A}_{\alpha}, \tilde{B}_{\alpha})$.
	Then, for any $k \in [n]$, the $(m+k)$th equation implies that
	\begin{align*}
		\alpha \otimes x_k \oplus \bigoplus_{i \neq k} x_i =\alpha \otimes x_k.
	\end{align*}
	This implies that $k \in M^{m+k}(\bm{x})$; hence, $\bigcup_{i=1}^{m+n} M^i(\bm{x})=[n]$.
	Proposition \ref{4-1} implies that $\bm{x}$ is a stable solution of 
	$\tilde{A}_{\alpha} \otimes \bm{x} = \tilde{B}_{\alpha} \otimes \bm{x}$.
\endproof

\subsection{Min-plus linear closure of the solution set}
\begin{figure}[htbp]
	\begin{center}
	\includegraphics[width=100mm]{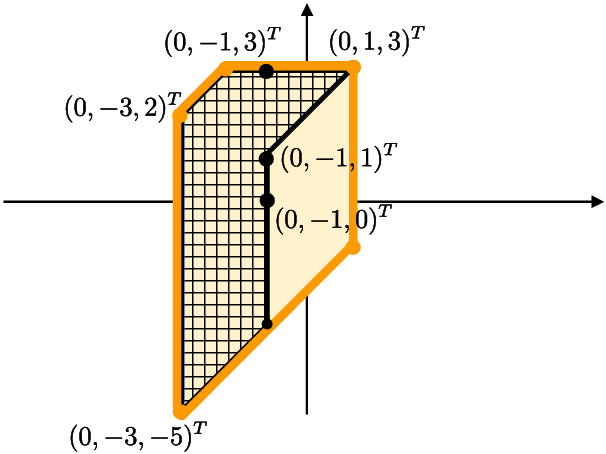}
	\caption{The solution set, $S(A,B)$, projected onto $x_{1}=0$ (lattice pattern) 
	and its min-plus linear closure (shaded).} \label{exmfig2}
	\end{center}
\end{figure}

In this section, we characterize the max-plus subspace $S(A,B)$ from a min-plus algebraic perspective.
Based on Proposition \ref{stableprop}, we consider the system $\tilde{A}_{\alpha} \otimes \bm{x} = \tilde{B}_{\alpha} \otimes \bm{x}$
instead of the system \eqref{hom-sys}, and assume that all solutions are stable, i.e., $S(A, B) = \tilde{S}(A, B)$.
Further, we assume that $S(A,B)$ is projectively bounded.
Let $\varphi(\bm{x})$ be the solution of \eqref{hom-sys} obtained by applying the alternating method to $\bm{x}(0) := \bm{x}$.
We use the notation $C = -(A \oplus B) \in \mathbb{R}_{\min}^{m \times n}$.
The $i$th row of $C$ is denoted by $C_{i}$.
\begin{proposition}\label{4-2}
	For any matrix $A,B\in\mathbb{R}_{\max}^{m\times n}$, we have
	\begin{align*}
		S(A,B) \subset \langle \varphi(C^T_1), \varphi(C^T_2), \dots, \varphi(C^T_m) \rangle_{\min} \cup \{\bm{\varepsilon}\}.
	\end{align*}
\end{proposition}
\proof
	Consider any vector $\bm{x} \in S(A,B) \setminus \{\bm{\varepsilon}\}$.
	Because $\bm{x}$ is a stable solution, the equality \eqref{stable} implies
	\begin{align*}
		\bm{x}=\bigoplus_{i=1}^m\!\text{\large ${}'$ } t_i \otimes' C^T_i,
	\end{align*}
	where $t_{i}$ denotes the $i$th entry of $A \otimes \bm{x}$ (or, equivalently, the $i$th entry of $B \otimes \bm{x}$).
	Based on the monotonicity of the operator $\varphi$, we have $\bm{x} = \varphi(\bm{x}) \leq \varphi(t_{i} \otimes' C_i^T)$.
	Hence, we have
	\begin{align*}
		\bm{x} \leq \bigoplus_{i=1}^m\!\text{\large ${}'$ } \varphi(t_i \otimes' C^T_i)
			= \bigoplus_{i=1}^m\!\text{\large ${}'$ } t_i \otimes' \varphi(C^T_i). 
	\end{align*}
	To demonstrate that $\varphi(C_{i}^{T}) \leq C_{i}^{T}$, it is sufficient to derive $\varphi_{0}(C_{i}^{T}) \leq C_{i}^{T}$
	because $\varphi(C_{i}^{T}) \leq \varphi_{0}(C_{i}^{T})$ holds by Lemma \ref{1-1}.
	We first note that $(A_i \otimes C_i^T) \oplus' (B_i \otimes C^T_i)=0$ as $C_i=-(A_i\oplus B_i)$.
	Then, we have
	\begin{align*}
		\varphi_{0}(C_{i}^{T}) &= C^T \otimes' ((A \otimes C_i^T) \oplus' (B \otimes C^T_i))\\
		&= \bigoplus_{k=1}^m\!\text{\large ${}'$ } ((A_{k} \otimes C_i^T) \oplus' (B_{k} \otimes C^T_i)) \otimes' C_k^T\\
		&= C_{i}^{T} \oplus' 
			\bigoplus_{k\neq i}\!\text{\large ${}'$ } ((A_{k} \otimes C_i^T) \oplus' (B_{k} \otimes C^T_i)) \otimes' C_k^T \\
		&\leq C_i^T.
	\end{align*}
	Thus, we have $\varphi(C_{i}^{T}) \leq C_{i}^{T}$; hence,
	\begin{align*}
		\bigoplus_{i=1}^m\!\text{\large ${}'$ } t_i \otimes' \varphi(C^T_i) 
			\leq \bigoplus_{i=1}^m\!\text{\large ${}'$ } t_i \otimes'C^T_i = \bm{x}.
	\end{align*}
	This proves the equality
	\begin{align*}
		\bm{x} = \bigoplus_{i=1}^m\!\text{\large ${}'$ } t_i \otimes' \varphi(C^T_i),
	\end{align*}
	which implies that $\bm{x} \in \langle \varphi(C^T_1), \varphi(C^T_2), \dots, \varphi(C^T_m) \rangle_{\min}$.
	Hence, we conclude that 
	$S(A,B) \subset \langle \varphi(C^T_1), \varphi(C^T_2), \dots, \varphi(C^T_m) \rangle_{\min} \cup \{\bm{\varepsilon}\}$,
\endproof
For a projectively bounded max-plus subspace $V \subset \mathbb{R}_{\max}^{n}$, 
the smallest min-plus subspace of $\mathbb{R}_{\min}^{n}$ containing $V \setminus \{\bm{\varepsilon}\}$ 
is defined to be the min-plus linear closure of $V$.
A projectively bounded max-plus subspace $V \subset \mathbb{R}_{\max}^{n}$ is called min-plus linear 
if $(V \setminus \{\bm{\varepsilon}\}) \cup \{\bm{\varepsilon}'\}$ is a min-plus subspace of $\mathbb{R}_{\min}^{n}$.
The following corollary asserts that the min-plus linear closure of the solution set $S(A,B)$ can be obtained via the alternating method.
\begin{corollary}\label{4-3}
	For any matrix $A,B\in\mathbb{R}_{\max}^{m\times n}$, the min-plus linear closure of $S(A,B)\setminus\{\bm{\varepsilon}\}$ is
	\begin{align*}
		\langle\varphi(C^T_1),\varphi(C^T_2),\ldots,\varphi(C^T_m)\rangle_{\min}.
	\end{align*}
	In particular, if $S(A,B)$ is min-plus linear, then
	\begin{align*}
		S(A,B) \setminus \{\bm{\varepsilon}\} 
		= \langle\varphi(C^T_1),\varphi(C^T_2),\ldots,\varphi(C^T_m)\rangle_{\min} \setminus \{\bm{\varepsilon}'\}.
	\end{align*}
\end{corollary}
\begin{example} \label{exm2}
	We consider the following matrices
	\begin{align*}
		A = \begin{pmatrix} 0  & 1 & -1 \\ 0 & -5 & -5 \\ 0 & 4 & 6 \\ 0 & 3 & -2 
			\\ \alpha & 0 & 0 \\ 0 & \alpha & 0 \\ 0 & 0 & \alpha \end{pmatrix}, \quad
		B = \begin{pmatrix} 0 & -1 & -1 \\ 0 & -4 & -3 \\ -1 & 1 & 6 \\ -1 & 3 & -3 
			\\ \alpha & -1 & -1 \\ -1 & \alpha & -1 \\ -1 & -1 & \alpha \end{pmatrix},
	\end{align*} 
	obtained from the matrices in Example \ref{exm1} via extension \eqref{extmat}.
	We set $\alpha = 13$.
	Then, by setting
	\begin{align*}
		C_{1} = (0,-1,1),& \quad C_{2} = (0,4,3), \quad C_{3} = (0,-4,-6), \quad C_{4} = (0,-3,2), \\
		C_{5} &= (\alpha,0,0), \quad C_{6} = (0,\alpha,0), \quad C_{7} = (0,0,\alpha),
	\end{align*}
	we have
	\begin{align*}
		&\varphi(C_{1}) = (0,-1,1)^{T}, \quad \varphi(C_{2}) = (0,1,3)^{T}, \\
		&\varphi(C_{3}) \sim (0,-3,-5)^{T}, \quad \varphi(C_{4}) = (0,-3,2)^{T}, \\
		\varphi(C_{5}) \sim (&0,-3,-5)^{T}, \quad \varphi(C_{6}) \sim (0,-1,0)^{T}, \quad \varphi(C_{7}) \sim (0,-1,3)^{T},
	\end{align*}
	where scalar multiples are modified for $C_{3}, C_{5}, C_{6}$, and $C_{7}$ so that the first entries are $0$.
	The min-plus subspace generated by these seven vectors is illustrated by the shaded region in Figure \ref{exmfig2}.
\end{example}
\section{A criterion for the solution set to be min-plus linear}
\begin{figure}[htbp]
	\begin{center}
	\includegraphics[width=60mm]{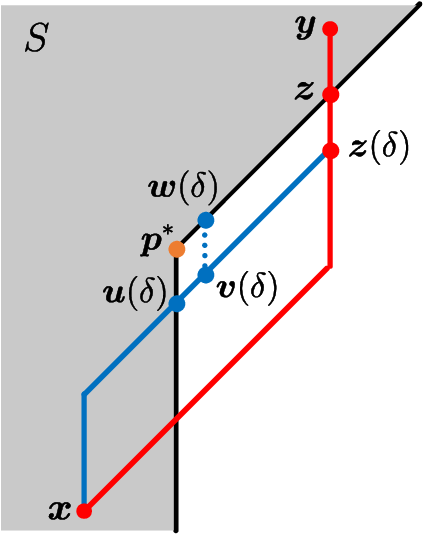}
	\caption{Proof of Lemma \ref{loc-glob}.} \label{proof52}
	\end{center}
\end{figure}

By Corollary \ref{4-3}, the solution set $S(A,B)$ can be computed via the alternating method if it is min-plus linear.
In this section, we present a sufficient condition for the solution set $S(A,B)$ to be min-plus linear.
First, we introduce the concept of local min-plus convex sets.
We define the distance between $\bm{x} = (x_{1},x_{2},\dots,x_{n})^{T} \in \mathbb{R}^{n}$
and $\bm{y} = (y_{1},y_{2},\dots,y_{n})^{T} \in \mathbb{R}^{n}$ as follows:
\begin{align*}
	d(\bm{x}, \bm{y}) = \max_{i \in [n]} |x_{i} - y_{i}|.
\end{align*}
The $r$-neighborhood of $\bm{x} = (x_{1},x_{2},\dots,x_{n})^{T} \in \mathbb{R}^{n}$ is defined as 
\begin{align*}
	\mathcal{N}(\bm{x}, r) = \{ \bm{y} \in \mathbb{R}^{n} \ |\ d(\bm{x}, \bm{y}) < r \}
\end{align*}
for $r > 0$.
\begin{definition} 
	Let $S \subset \mathbb{R}^{n} \cup \{\bm{\varepsilon}\}$ be a subset closed under scalar multiplication.
	Then, $S$ is said to be locally min-plus convex at $\bm{x} \in S \setminus \{\bm{\varepsilon}\}$ if 
	there exists $r>0$ such that $\bm{y} \oplus' \bm{z} \in S$ holds for any $\bm{y}, \bm{z}\in \mathcal{N}(\bm{x},r) \cap S$.
\end{definition}
In fact, the min-plus local convexity at every point of the max-plus subspace is equivalent to global min-plus linearity.
\begin{proposition} \label{loc-glob}
	Let $S \subset \mathbb{R}_{\max}^{n}$ be a projectively bounded max-plus subspace.
	Then, $S$ is min-plus linear if and only if it is locally min-plus convex for any $\bm{x} \in S \setminus \{\bm{\varepsilon}\}$.
\end{proposition}
\proof
	The ``only if'' part is trivial.
	We need to prove the ``if'' part of the proposition.
	Suppose that $S$ is not min-plus linear.
	Then, there exist two vectors $\bm{x}, \bm{y} \in S \setminus \{\bm{\varepsilon}\}$ 
	such that $\bm{x} \oplus' \bm{y} \not\in S$.
	We will find the point $\bm{p}^{*} \in S$ where $S$ is not locally min-plus convex, as shown in Figure \ref{proof52}.
	Consider the minimum positive number $\alpha$ such that
	\begin{align*}
		\bm{z} := \alpha \otimes' \bm{x} \oplus' \bm{y} \in S.
	\end{align*}
	By minimality of $\alpha$, we have
	\begin{align*}
		\bm{z}(\delta) := (\alpha-\delta) \otimes' \bm{x} \oplus' \bm{y} \not\in S
	\end{align*}
	when $0 < \delta < \alpha$.
	Further, let $\beta(\delta)$ be the maximum number such that
	\begin{align*}
		\bm{x} \oplus \beta(\delta) \otimes \bm{z}(\delta) \in S.
	\end{align*}
	We set $\bm{u}(\delta) = \delta \otimes (\bm{x} \oplus \beta(\delta) \otimes \bm{z}(\delta))$.
	By maximality of $\beta(\delta)$, we have
	\begin{align*}
		\bm{v}(\delta) := \bm{x} \oplus (\beta(\delta) + \delta) \otimes \bm{z} (\delta) \not\in S.
	\end{align*}
	As $S$ is a max-plus subspace, we have
	\begin{align*}
		\bm{w}(\delta) := \bm{x} \oplus (\beta(\delta) + \delta) \otimes \bm{z} \in S.
	\end{align*}
	Now, we prove that $\bm{u}(\delta), \bm{v}(\delta)$ and $\bm{w}(\delta)$ are contained in 
	the neighborhood of a certain fixed point $\bm{p}^{*}$ and satisfy $\bm{v}(\delta) = \bm{u}(\delta) \oplus' \bm{w}(\delta)$.
	First, we determine $\bm{p}^{*}$ by demonstrating that $\lim_{\delta \to +0} \beta(\delta)$ exists.
	\par
\begin{lemma} \label{betalem}
	If $0 < \delta' < \delta$, then $\beta(\delta) \leq \beta(\delta')$.
	Further, when $0 < \delta < 1$, $\beta(\delta)$ is bounded above.
\end{lemma}
\proof
	Let $\bm{x} = (x_{1},x_{2},\dots,x_{n})^{T}$ and $\bm{z} = (z_{1},z_{2},\dots,z_{n})^{T}$.
	Firstly, there exist $\delta > 0$ and $I \subset [n]$ such that
	\begin{align*}
		z_{i}(\delta') = \begin{cases} z_{i} & i \in I, \\ z_{i} - \delta' & i \not\in I, \end{cases}
	\end{align*}
	for any $0 < \delta' \leq \delta$.
	Then, we have 
	\begin{align*}
		\bm{z}(\delta') = \bm{z}(\delta) \oplus (-\delta') \otimes \bm{z}.
	\end{align*}
	Because $\bm{x} \oplus \beta(\delta) \otimes \bm{z}(\delta) \in S$, we have
	\begin{align*}
		\bm{x} \oplus \beta(\delta) \otimes \bm{z}(\delta') 
			= (\bm{x} \oplus \beta(\delta) \otimes \bm{z}(\delta)) \oplus (\beta(\delta) - \delta') \otimes \bm{z} \in S.
	\end{align*}
	By maximality of $\beta(\delta')$, we have $\beta(\delta) \leq \beta(\delta')$.
	Next, let $\xi = \max_{i} x_{i}$ and $\zeta = \min_{i} z_{i}$.
	Because $\bm{z}(\delta) \geq (-1) \otimes \bm{z}$ for $\delta < 1$, we have:
	\begin{align*}
		(\xi - \zeta + 1) \otimes \bm{z}(\delta) \geq \xi \otimes (0,0,\dots,0)^{T} \geq \bm{x}.
	\end{align*}
	Hence, for any $\beta \geq \xi - \zeta + 1$, we have:
	\begin{align*}
		\bm{x} \oplus \beta \otimes \bm{z}(\delta) = \beta \otimes \bm{z}(\delta) \not\in S,
	\end{align*}
	which implies that $\beta(\delta) < \xi - \zeta + 1$,
\endproof
\proof[Proof of Proposition \ref{loc-glob} (continued)]
	By Lemma \ref{betalem}, let $\beta^{*} = \lim_{\delta \to +0} \beta(\delta)$.
	Then, the desired point $\bm{p}^{*}$ is given by
	\begin{align*}
		\bm{p}^{*} = \bm{x} \oplus \beta^{*} \otimes \bm{z}.
	\end{align*}
	Indeed, for any $r > 0$, we can verify that $\bm{u}(\delta), \bm{v}(\delta), \bm{w}(\delta) \in \mathcal{N}(\bm{p}^{*},r)$
	by taking sufficiently small $\delta > 0$ such that $\beta^{*} - \beta(\delta) < \frac{r}{2}$ and $\delta < \frac{r}{2}$.
	\par
	Next, we demonstrate that $\bm{v}(\delta) = \bm{u}(\delta) \oplus' \bm{w}(\delta)$.
	Because $\bm{z}(\delta) \leq (\alpha-\delta) \otimes' \bm{x}$, we have
	\begin{align*}
		\bm{x} \oplus (\delta-\alpha) \otimes \bm{z}(\delta) = \bm{x} \in S.
	\end{align*}
	By maximality of $\beta(\delta)$, we have $\alpha + \beta(\delta) \geq \delta$.
	We also note that
	\begin{align*}
		\bm{z}(\delta) = (\alpha-\delta) \otimes' \bm{x} \oplus' \bm{y}
			= (\alpha-\delta) \otimes' \bm{x} \oplus' \alpha \otimes' \bm{x} \oplus' \bm{y}
			=  (\alpha-\delta) \otimes' \bm{x} \oplus' \bm{z}.
	\end{align*}
	Then, we have
	\begin{align*}
		\bm{u}(\delta) &= \delta \otimes \bm{x} 
			\oplus (\beta(\delta) + \delta) \otimes ((\alpha-\delta) \otimes' \bm{x} \oplus' \bm{z}) \\
		&= (\delta \otimes \bm{x} \oplus (\alpha + \beta(\delta)) \otimes \bm{x}) 
			\oplus' (\delta \otimes \bm{x} \oplus (\beta(\delta) + \delta) \otimes \bm{z}) \\
		&= ((\alpha + \beta(\delta)) \otimes \bm{x}) \oplus' (\delta \otimes \bm{x} \oplus (\beta(\delta) + \delta) \otimes \bm{z}).
	\end{align*}
	Here, we use the property $\bm{a} \oplus (\bm{b} \oplus' \bm{c}) = (\bm{a} \oplus \bm{b}) \oplus' (\bm{a} \oplus \bm{c})$.
	Thus, we have
	\begin{align*}
		& \bm{u}(\delta) \oplus' \bm{w}(\delta) \\
		=\, & ((\alpha + \beta(\delta)) \otimes \bm{x}) \oplus' (\delta \otimes \bm{x} \oplus (\beta(\delta) + \delta) \otimes \bm{z})
			\oplus' (\bm{x} \oplus (\beta(\delta) + \delta) \otimes \bm{z}) \\
		=\, & ((\alpha + \beta(\delta)) \otimes \bm{x}) \oplus' (\bm{x} \oplus (\beta(\delta) + \delta) \otimes \bm{z}).
	\end{align*}
	On the other hand, we can express $\bm{v}(\delta)$ as follows
	\begin{align*}
		\bm{v}(\delta) &= \bm{x} \oplus (\beta(\delta) + \delta) \otimes ((\alpha-\delta) \otimes' \bm{x} \oplus' \bm{z}) \\
		&= (\bm{x} \oplus (\alpha + \beta(\delta)) \otimes \bm{x}) \oplus' (\bm{x} \oplus (\beta(\delta)+\delta) \otimes \bm{z}) \\
		&= ((\alpha + \beta(\delta)) \otimes \bm{x}) \oplus' (\bm{x} \oplus (\beta(\delta)+\delta) \otimes \bm{z}),
	\end{align*}
	This proves that $\bm{v}(\delta) = \bm{u}(\delta) \oplus' \bm{w}(\delta)$.
	Hence, we conclude that $S(A,B)$ is not min-plus convex at $\bm{p}^{*} \in S$.
\endproof

Next, we discuss the local min-plus convexity of the solution set $S(A, B)$.
\begin{definition}\label{4-4}
	Let $\bm{a}=(a_1,a_2,\dots,a_n),\bm{b}=(b_1,b_2,\dots,b_n)\in\mathbb{R}^{1\times n}_{\max}$.
	We define the subset $R(\bm{a}, \bm{b}) \subset \mathbb{R}^{n}$
	consisting of solutions $\bm{x} = (x_{1},x_{2},\dots,x_{n})^{T}$
	 of a single equation $\bm{a} \otimes \bm{x} = \bm{b} \otimes \bm{x}$ satisfying one of the following conditions:
	\begin{enumerate}
	\item
	There exist $i,j\in[n],\ i\neq j,$ such that
	\begin{align*}
	\bm{a}\otimes\bm{x} &=a_i\otimes x_i>a_k\otimes x_k\quad \text{ for } k\neq i, \text{ and}\\
	\bm{b}\otimes\bm{x} &=b_j\otimes x_j>b_k\otimes x_k\quad \text{ for } k\neq j.
	\end{align*}
	\item
	There exists $i \in [n]$ such that
	\begin{align*}
	\bm{a}\otimes\bm{x} &=a_i\otimes x_i, \\
	\bm{b}\otimes\bm{x} &=b_i\otimes x_i>b_k\otimes x_k\quad \text{ for } k\neq i.
	\end{align*}
	\item
	For any $i \in [n]$,
	\begin{align*}
	\bm{a}\otimes\bm{x}=a_i\otimes x_i \quad\Longleftrightarrow\quad \bm{b}\otimes\bm{x}=b_i\otimes x_i
	\end{align*}
	holds.
	\end{enumerate}
\end{definition}
\begin{proposition}\label{kyokusho}
	Let $\bm{a}=(a_1,a_2,\dots,a_n),\bm{b}=(b_1,b_2,\dots,b_n)\in\mathbb{R}^{1\times n}_{\max} \setminus \{\bm{\varepsilon}\}$.
	For any vector $\bm{x} = (x_{1},x_{2},\dots,x_{n})^{T} \in S(\bm{a},\bm{b}) \cap \mathbb{R}^{n}$, (i) and (ii) are equivalent:
	\begin{enumerate}
	\item[(i)]
		$S(\bm{a},\bm{b})$ is locally min-plus convex at $\bm{x}$.
	\item[(ii)]
		$\bm{x} \in R(\bm{a},\bm{b})$.
	\end{enumerate}
\end{proposition}
\proof
(ii) $\Longrightarrow$ (i): 
Take any vector $\bm{x}\in R(\bm{a},\bm{b})$ and define
\begin{align*}
	K(\bm{a},\bm{x}) &=\{i \in [n] \mid \bm{a} \otimes \bm{x} = a_i \otimes x_i\},\\
	K(\bm{b},\bm{x}) &=\{i \in [n] \mid \bm{b} \otimes \bm{x} = b_i \otimes x_i\}.
\end{align*}
Then, for sufficiently small $r > 0$ and $\bm{y}, \bm{z} \in \mathcal{N}(\bm{x},r) \cap S(\bm{a},\bm{b})$, we have
\begin{align*}
	a_i\otimes y_i >a_k\otimes y_k, \quad a_i\otimes z_i >a_k\otimes z_k
\end{align*}
for $i\in K(\bm{a},\bm{x}), k\notin K(\bm{a},\bm{x})$, and
\begin{align*}
	b_i\otimes y_i >b_k\otimes y_k, \quad b_i\otimes z_i &>a_k\otimes z_k
\end{align*}
for $i\in K(\bm{b},\bm{x}), k\notin K(\bm{b},\bm{x})$.
This yields
\begin{align*}
	K(\bm{a},\bm{x}) \supset K(\bm{a},\bm{y}), K(\bm{a},\bm{z}), \quad
	K(\bm{b},\bm{x}) \supset K(\bm{b},\bm{y}), K(\bm{b},\bm{z}).
\end{align*}
\begin{itemize}
\item Suppose $\bm{x}$ satisfies condition 1 of Definition \ref{4-4} for $K(\bm{a},\bm{x})=\{i\}$ and $K(\bm{b},\bm{x})=\{j\}$.
Then, we have
\begin{align*}
	a_i \otimes (y_i \oplus' z_i) &> a_k \otimes (y_k \oplus' z_k) \quad \text{ for } k\neq i, \\
	b_j \otimes (y_j \oplus' z_j) &> a_k \otimes (y_k \oplus' z_k) \quad \text{ for } k\neq j.
\end{align*}
Therefore, we compute
\begin{align*}
	\bm{a} \otimes (\bm{y} \oplus' \bm{z})
	&= a_i \otimes (y_i \oplus' z_i)\\
	&= (a_i \otimes y_i) \oplus' (a_i \otimes z_i)\\
	&= (b_j \otimes y_j) \oplus' (b_j \otimes z_j)\\
	&= b_j \otimes (y_j \oplus' z_j)\\
	&= \bm{b} \otimes(\bm{y} \oplus' \bm{z})
\end{align*}
This proves that $\bm{y} \oplus' \bm{z} \in S(\bm{a},\bm{b})$.
\item Suppose $\bm{x}$ satisfies condition 2 of Definition \ref{4-4} for $K(\bm{b},\bm{x})=\{i\}$,
which yields $a_{i} = b_{i}$.
Then, we have
\begin{align*}
	a_i\otimes y_i=b_i\otimes y_i&>b_k\otimes y_k \quad \text{ for } k\neq i, \\
	a_i\otimes z_i=b_i\otimes z_i&>b_k\otimes z_k \quad \text{ for } k\neq i, 
\end{align*}
and hence
\begin{align*}
	a_i\otimes y_i&\geq a_k \otimes y_k \quad \text{ for } k\neq i,  \\
	a_i\otimes z_i&\geq a_k \otimes z_k \quad \text{ for } k\neq i.
\end{align*}
Therefore, we have
\begin{align*}
	\bm{a} \otimes (\bm{y} \oplus'\bm{z}) = a_i \otimes (y_i \oplus' z_i)
	= b_i \otimes (y_i \oplus' z_i)
	= \bm{b} \otimes (\bm{y} \oplus' \bm{z}).
\end{align*}
Thus, we see that $\bm{y}\oplus'\bm{z}\in S(\bm{a},\bm{b})$.
\item Suppose $\bm{x}$ satisfies condition 3 of Definition \ref{4-4}.
Then, we have
\begin{align*}
	K(\bm{a},\bm{y})=K(\bm{b},\bm{y}), \quad K(\bm{a},\bm{z})=K(\bm{b},\bm{z})
\end{align*}
since $a_{i} = b_{i}$ for all $i \in K(\bm{a},\bm{x})=K(\bm{b},\bm{x})$.
Therefore, we have
\begin{align*}
	\bm{a}\otimes(\bm{y}\oplus'\bm{z})&=a_i\otimes(y_i\oplus'z_i)\\
	&=(a_i\otimes y_i)\oplus'(a_i\otimes z_i)\\
	&=(b_i\otimes y_i)\oplus'(b_i\otimes z_i)\\
	&=b_i\otimes (y_i\oplus'z_i)\\
	&=\bm{b}\otimes(\bm{y}\oplus'\bm{z}).
\end{align*}
Thus, we see that $\bm{y}\oplus'\bm{z}\in S(\bm{a},\bm{b})$.
\end{itemize}
Based on these three cases, we conclude that $S(\bm{a}, \bm{b})$ is locally min-plus convex at $\bm{x}$.
\par
(i) $\Longrightarrow$ (ii): 
Consider any vector $\bm{x} \in \mathbb{R}^{n}$ such that $\bm{x} \not\in R(\bm{a}, \bm{b})$.
Then, we have the following four cases:
\begin{enumerate}
\item[(a)] When $|K(\bm{a},\bm{x})|=|K(\bm{b},\bm{x})|=1$, 
$K(\bm{a},\bm{x}) = K(\bm{b},\bm{x})$ implies condition 3 of Definition \ref{4-4} and 
$K(\bm{a},\bm{x})\neq K(\bm{b},\bm{x})$ implies condition 1 of Definition \ref{4-4}. 
Hence, this case cannot occur if  $\bm{x} \not\in R(\bm{a}, \bm{b})$.
\item[(b)] When $|K(\bm{a},\bm{x})| = 1$ and $|K(\bm{b},\bm{x})| \neq 1$,
$K(\bm{a},\bm{x}) \subset K(\bm{b},\bm{x})$ implies condition 2 in Definition \ref{4-4}. 
Hence, we have $K(\bm{a},\bm{x}) \not\subset K(\bm{b},\bm{x})$.
\item[(c)] When $|K(\bm{a},\bm{x})| \neq 1$ and $|K(\bm{b},\bm{x})| = 1$,
$K(\bm{b},\bm{x}) \not\subset K(\bm{a},\bm{x})$ as in case (b).
\item[(d)] When $|K(\bm{a},\bm{x})| \neq 1$ and $|K(\bm{b},\bm{x})| \neq 1$,
$K(\bm{a},\bm{x})=K(\bm{b},\bm{x})$ implies condition 3 of Definition \ref{4-4}. 
Hence, we have $K(\bm{a},\bm{x})\neq K(\bm{b},\bm{x})$.
\end{enumerate}
From the aforementioned observations, there exist three distinct indices, $i_1,i_2,j$, such that:
\begin{align*}
	\left\{
	\begin{array}{l}
	j\in K(\bm{b},\bm{x})\setminus K(\bm{a},\bm{x}),\\
	i_1,i_2\in K(\bm{a},\bm{x}),
	\end{array}
	\right.
\end{align*}
or
\begin{align*}
	\left\{
	\begin{array}{l}
	j\in K(\bm{a},\bm{x})\setminus K(\bm{b},\bm{x}),\\
	i_1,i_2\in K(\bm{b},\bm{x}).
	\end{array}
	\right.
\end{align*}
In the former case, 
defining the two vectors $\bm{y}_{1}, \bm{y}_{2}$ as 
\begin{align*}
	\bm{y}_1=(x_1,x_2,\dots,x_{i_1-1},x_{i_1}+r, x_{i_1+1},\dots,x_{j-1},x_{j}+r, x_{j+1},\dots,x_n)^T, \\
	\bm{y}_2=(x_1,x_2,\dots,x_{i_2-1},x_{i_2}+r, x_{i_2+1},\dots,x_{j-1},x_{j}+r, x_{j+1},\dots,x_n)^T,
\end{align*}
for a sufficiently small $r > 0$, we have
\begin{align*}
	\bm{a} \otimes \bm{y}_1 &= \bm{a} \otimes \bm{x} + r = \bm{b} \otimes \bm{x} + r = \bm{b} \otimes \bm{y}_1, \\
	\bm{a} \otimes \bm{y}_2 &= \bm{a} \otimes \bm{x} + r = \bm{b} \otimes \bm{x} + r = \bm{b} \otimes \bm{y}_2.
\end{align*}
Hence, $\bm{y}_1,\bm{y}_2\in R(\bm{x},r) \cap S(\bm{a},\bm{b})$.
However, the vector $\bm{z} = \bm{y}_1\oplus' \bm{y}_2$ satisfies
$\bm{a} \otimes \bm{z} = \bm{a} \otimes \bm{x}$ and $\bm{b} \otimes \bm{z} = \bm{b} \otimes \bm{x} + r$,
which yields $\bm{a} \otimes \bm{z} \neq \bm{b} \otimes \bm{z}$.
Because $r > 0$ can be taken to be arbitrarily small, $S(A,B)$ cannot be locally min-plus convex at $\bm{x}$.
\endproof
To verify the min-plus linearity of the solution set $S(A,B)$, the local min-plus convexity at only some special points needs to be verified.
First, we present the following technical lemma.
\begin{lemma}\label{junnbi}
	If $\bm{y} \in S(A,B)$ satisfies $\bm{y}\leq\bm{x}$ for $\bm{x}\in\mathbb{R}^n$, then $\bm{y} \leq \varphi(\bm{x})$.
\end{lemma}
\proof
	We denote by $\bm{x}(r)$ the vector in the $r$th iteration of the alternating method applied to $\bm{x}(0) := \bm{x}$.
	It is sufficient to demonstrate that $\bm{y}\leq\bm{x}(r)$ implies $\bm{y}\leq\bm{x}(r+1)$.
	First, we assume  $\bm{y}\leq\bm{x}(r)$ and $\bm{z}=A\otimes\bm{y}=B\otimes\bm{y}$.
	Then, we have
	\begin{align*}
		\bm{z} = (A \otimes \bm{y}) \oplus' (B \otimes \bm{y}) \leq (A \otimes \bm{x}(r)) \oplus' (B \otimes \bm{x}(r)),
	\end{align*}
	which yields
	\begin{align*}
		C^T\otimes'\bm{z}\leq C^T\otimes'((A\otimes\bm{x}(r))\oplus'(B\otimes\bm{x}(r))). 
	\end{align*}
	As $(-A^T) \otimes' (A \otimes \bm{y}) \geq \bm{y}$ and $(-B^T) \otimes' (B \otimes \bm{y}) \geq \bm{y}$,
	 we have
	\begin{align*}
		C^T\otimes'\bm{z} = (-A^T)\otimes'(A\otimes\bm{y})\oplus'(-B^T)\otimes'(B\otimes\bm{y}) \geq \bm{y}. 
	\end{align*}
	This proves that
	\begin{align*}
		\bm{y} \leq C^T\otimes'((A\otimes\bm{x}(r))\oplus'(B\otimes\bm{x}(r))) = \bm{x}(r+1).
	\end{align*}
	We obtain $\bm{y} \leq \varphi(\bm{x})$ via induction.
\endproof
Now, we present the main results.
\begin{theorem} \label{hantei}
	If $\varphi(C_i^T) \in R(A_i,B_i)$ holds for all $i\in[m]$, then $S(A,B)$ is min-plus linear.
\end{theorem}
\proof
	Let us assume that $S(A,B)$ is not min-plus linear. 
	Then, by Proposition \ref{loc-glob}, 
	there exists a vector $\bm{x}^{\star} = (x_{1}^{\star}, x_{2}^{\star}, \dots, x_{n}^{\star})^{T} \in S(A,B)$ 
	such that $S(A,B)$ is not locally min-plus convex at $\bm{x}^{\star}$.
	In this case, $S(A_{i}, B_{i})$ is not locally min-plus convex at $\bm{x}^{\star}$ for some $i \in [m]$,
	which implies that $\bm{x}^{\star} \notin R(A_i,B_i)$ by Proposition \ref{kyokusho}. 
	It may be assumed that $\bm{x}^{\star}$ satisfies $\bm{x}^{\star}\leq C_i^T$ and $x^{\star}_k=c_{ik}$ for some $k \in [n]$ 
	because $R(A_{i}, B_{i})$ is closed under scalar multiplication.
	Because $(A_{i} \oplus B_{i}) \otimes C_i^T = 0 = (a_{ik} \oplus b_{ik}) \otimes c_{ik}$, we have
	\begin{align*}
		A_i\otimes\bm{x}^{\star} = B_i\otimes\bm{x}^{\star} = (A_i \oplus B_{i}) \otimes C_i^T.
	\end{align*}
	Further, because $\bm{x}^{\star} \leq \varphi(C_i^T) \leq C_i^T$ by the proof of Proposition \ref{4-2} and Lemma \ref{junnbi},
	we have
	\begin{align*}
		A_i \otimes \bm{x}^{\star} = B_i \otimes \bm{x}^{\star} = A_i \otimes \varphi(C_i^T) = B_i \otimes \varphi(C_i^T)
		= (A_i \oplus B_{i}) \otimes C_i^T.
	\end{align*}
	Therefore, we obtain
	\begin{align*}
		K(A_i,\bm{x}^{\star}) \subset K(A_i, \varphi(C^T_i)),\quad K(B_i, \bm{x}^{\star}) \subset K(B_i, \varphi(C^T_i)).
	\end{align*}
	Hence, we can easily verify that $\bm{x}^{\star} \notin R(A_i,B_i)$ leads to $\varphi(C^T_i) \notin R(A_i,B_i)$.
	This contradicts the assumption that $\varphi(C^T_i) \in R(A_i,B_i)$ holds for all $i \in [m]$.
	This concludes the proof of this theorem.
\endproof
We note that Theorem \ref{hantei} provides a sufficient, but not necessary, condition for $S(A,B)$ to be min-plus linear.
This is due to the fact that $S(A,B)$ is locally min-plus convex at $\bm{x}$ if $S(A_{i},B_{i})$ is so for all $i \in [m]$,
but the converse is not true.
However, we can use Theorem \ref{hantei} as a rather strict criterion.
\par
\begin{example}
	Consider the linear system \eqref{hom-sys} for
	\begin{align*}
		A = \begin{pmatrix} 0  & 1 & -1 \\ 0 & -5 & -5 \\ 0 & 4 & 6 \\ 0 & 3 & -2 \end{pmatrix}, \quad
		B = \begin{pmatrix} 0 & -1 & -1 \\ 0 & -4 & 0 \\ -1 & 1 & 6 \\ -1 & 3 & -3 \end{pmatrix}.
	\end{align*} 
	By applying the alternating method, we obtain
	\begin{align*}
		&\varphi(C_{1}) = (0,-1,0)^{T}, \quad \varphi(C_{2}) = (0,-1,0)^{T}, \\
		&\varphi(C_{3}) = (0,-3,-5)^{T}, \quad \varphi(C_{4}) = (0,-3,0)^{T}.
	\end{align*}
	Since 
	\begin{align*}
		&K(A_{1},\varphi(C_{1})) = \{ 1,2 \}, \quad K(B_{1},\varphi(C_{1})) = \{ 1 \}, \\
		&K(A_{2},\varphi(C_{2})) = \{ 1 \}, \quad K(B_{2},\varphi(C_{2})) = \{ 1,3 \}, \\
		&K(A_{3},\varphi(C_{3})) = \{ 2,3 \}, \quad K(B_{3},\varphi(C_{3})) = \{ 3 \}, \\
		&K(A_{4},\varphi(C_{4})) = \{ 1,2 \}, \quad K(B_{3},\varphi(C_{3})) = \{ 2 \}, 
	\end{align*}
	the condition in Theorem \ref{hantei} is satisfied.
	Hence, $S(A,B)$ is min-plus linear and 
	\begin{align*}
		S(A,B) \setminus \{\bm{\varepsilon}\} = 
			\langle\varphi(C^T_1),\varphi(C^T_2),\varphi(C^T_3),\varphi(C^T_{4})\rangle_{\min} \setminus \{\bm{\varepsilon}'\}.
	\end{align*}
\end{example}
Summarizing the above results, we have the following algorithm that characterizes the solution set $S(A, B)$ of 
the max-plus linear system \eqref{hom-sys} for $A, B \in \mathbb{R}_{\max}^{m \times n}$.
It is not required for $S(A,B)$ to be projectively bounded nor stable solutions.
\begin{algorithm} \label{mainalg}
	\begin{enumerate}
	\item Let 
	\begin{align*}
		\tilde{A} =\begin{pmatrix}
		A\\
		D_{\alpha,0}
		\end{pmatrix}, \quad
		\tilde{B} =\begin{pmatrix}
		B\\
		D_{\alpha,-1}
		\end{pmatrix}.
	\end{align*}
	for $\alpha = \max\{ |c_{ij} - c_{ik}| \mid i \in [m],\ j ,k \in [n],\ c_{ki}, c_{kj} \in \mathbb{R} \} + 1$, 
	where $C = (c_{ij}) = -(A \oplus B)$.
	\item Replace all entries $\varepsilon$ in $\tilde{A}$ and $\tilde{B}$ with 
	\begin{align*}
	 	\beta = \min( \min\{ a_{ij} \mid a_{ij} \in \mathbb{R} \}, \min\{ b_{ij} \mid b_{ij} \in \mathbb{R} \})  - \alpha.
	\end{align*}
	\item Let $\tilde{C}_{i}$ be the $i$th row of $-(\tilde{A} \oplus \tilde{B})$.
	Compute $\varphi(\tilde{C}_{i}^{T})$ for all $i \in [m+n]$ via the alternating method.
	\item If $\left| [\varphi(\tilde{C}_{i}^{T})]_{j} - [\varphi(\tilde{C}_{i}^{T})]_{k} \right| < \alpha$ for all $i \in [m+n]$ and $j,k \in [n]$,
	then $S(A,B) = S(\tilde{A},\tilde{B})$ and it is projectively bounded; 
	otherwise $S(A,B)$ is not projectively bounded and it is approximated by the projectively bounded subspace $S(\tilde{A},\tilde{B})$.
	\item If $\varphi(\tilde{C}_{i}^{T}) \in R(A_{i}, B_{i})$ for all $i \in [m+n]$, 
	then $S(\tilde{A},\tilde{B})$ is min-plus linear.
	\item The min-plus linear closure of $S(\tilde{A},\tilde{B})$ is 
	\begin{align*}
		\langle\varphi(C^T_1),\varphi(C^T_2),\ldots,\varphi(C^T_{m+n})\rangle_{\min}.
	\end{align*}
	In particular, if $S(A,B)$ is projectively bounded and min-plus linear, then
	\begin{align*}
		S(A,B) \setminus \{\bm{\varepsilon}\} = 
			\langle\varphi(C^T_1),\varphi(C^T_2),\ldots,\varphi(C^T_{m+n})\rangle_{\min} \setminus \{\bm{\varepsilon}'\}.
	\end{align*}	
	\end{enumerate}
\end{algorithm}
By Lemma \ref{Dsol}, the replacement in Step 2 does not affect the set $S(\tilde{A},\tilde{B})$.
Steps 4 and 5 operate correctly using Proposition \ref{Dinc} and Theorem \ref{hantei}, respectively.
Step 6 follows from Corollary \ref{4-3}.
The complexity of this algorithm is as follows:
\begin{proposition}
	If $A = (a_{ij}), B = (b_{ij}) \in \mathbb{Z}^{m \times n}$, then the computational complexity of \ref{mainalg} is $O(n (m+n)^3 K)$.
	where $K = \max( \max_{i,j} |a_{ij}|, \max_{i,j} |b_{ij}|)$.
\end{proposition}
\proof
	Because $\alpha \leq 2K + 1$, all entries of $\tilde{A}$ and $\tilde{B}$ are bounded by $\pm (3K+1)$ 
	after replacement in Step 2.
	Therefore, the complexity of the alternating method for $\tilde{A} \otimes \bm{x} = \tilde{B} \otimes \bm{x}$ is $O(n (m+n)^{2} K)$.
	Step 3 requires $(m+n)$ iterations of the alternating method.
	Verification of Steps 4 and 5 requires $O(n^{2}(m+n))$ and $O(n(m+n))$ computation time, respectively.
	Thus, the total computational complexity of the algorithm is $O(n (m+n)^3 K)$.
\endproof


\end{document}